\documentclass[12pt]{amsart}
\usepackage{amssymb, amsfonts}

\newcommand{\<}{\langle}
\renewcommand{\>}{\rangle}

\renewcommand{\d}{\delta}

\newcommand{\e}{\varepsilon}

\renewcommand{\l}{\lambda}

\newcommand{\var}{\varphi}

\newcommand{\cA}{{\mathcal A}}
\newcommand{\cB}{{\mathcal B}}
\newcommand{\cC}{{\mathcal C}}
\newcommand{\cD}{{\mathcal D}}

\newcommand{\cF}{{\mathcal F}}
\newcommand{\cH}{{\mathcal H}}

\newcommand{\cL}{{\mathcal L}}

\newcommand{\cU}{{\mathcal U}}
\newcommand{\cV}{{\mathcal V}}

\newcommand{\cAo}{{\mathcal A}'_1}

\newcommand{\cFo}{{\mathcal F}'_1}

\newcommand{\cAH}{{\mathcal A}^h}
\newcommand{\cBH}{{\mathcal B}^h}

\newcommand{\cAD}{({\mathcal A^h})'_1}

\newcommand{\bR}{{\mathbb R}}
\newcommand{\bC}{{\mathbb C}}

\newtheorem{theorem}{Theorem}[section]
\newtheorem{corollary}[theorem]{Corollary} 
\newtheorem{proposition}[theorem]{Proposition} 

\newtheorem{lemma}[theorem]{Lemma} 
 
\newtheorem*{thmnonum}{Theorem}

\theoremstyle{definition}   
\newtheorem{definition}[theorem]{Definition} 
\newtheorem{example}[theorem]{Example}

\begin{document}

\title[Best Approximation from C*-subalgebras]{Leibniz Seminorms and 
Best Approximation from C*-subalgebras}
\author{Marc A. Rieffel}
\address{Department of Mathematics \\
University of California \\ Berkeley, CA 94720-3840}
\email{rieffel@math.berkeley.edu}
\thanks{The research reported here was
supported in part by National Science Foundation grant DMS-0753228.}
\dedicatory{In celebration of the successful completion by
Richard V. Kadison of 85 circumnavigations of the sun }


\begin{abstract}
We show that if $\cB$ is a C*-subalgebra of a C*-algebra
$\cA$ such that $\cB$ contains a bounded approximate identity
for $\cA$, and if $L$ is the pull-back to $\cA$ of the 
quotient norm on $\cA/\cB$, then
$L$ is strongly Leibniz. In connection with this situation
we study certain aspects of best approximation of elements of
a unital C*-algebra by elements of a unital C*-subalgebra. 
\end{abstract}

\maketitle
\allowdisplaybreaks

\section{Introduction}

From my attempts to understand C*-metrics, as defined in \cite{R21}, I have
recently been trying to discover the mechanisms that can lead to (continuous) 
seminorms, $L$,
defined on a C*-algebra $\cA$, that are Leibniz, that is, satisfy the Leibniz 
inequality
\[
L(AC) \leq L(A)\|C\| +\|A\|L(C)
\]
for all $A, C \in \cA$. Also of much interest to me, because of their importance
in \cite{R17} and in potential non-commutative versions of the results in
\cite{R17}, are seminorms on a unital C*-algebra that are strongly-Leibniz,
that is, satisfy, in addition to the Leibniz inequality, the property that if $A \in \cA$
and if $A$ is invertible in $\cA$, then
\[
L(A^{-1}) \leq \|A^{-1}\|^2L(A)  ,
\]
and also $L(1_A) = 0$, again as defined in \cite{R21}. This last inequality seems
to have received virtually no attention in the mathematics literature.

Actually, for infinite-dimensional C*-algebras, the C*-metrics as defined in\cite{R21}
are discontinuous and only densely defined. But they are required
to be lower semi-continuous with
respect to the C*-norm, and in all of the examples that I know of one proves that they are
lower semi-continuous by showing that they are the supremum of an infinite
family of continuous strongly-Leibniz seminorms. This provides ample reason for
studying continuous strongly-Leibniz seminorms.

My investigations have led me to consider the situation in which $\cB$ is a C*-subalgebra 
of a C*-algebra $\cA$ and $L$ is the quotient norm on $\cA/\cB$ pulled back to $\cA$.
That is,
\[
L(A) = \inf \{\|A-B\|\  :  B \in \cB\}  
\]
for $A \in \cA$, so that $L(A)$ is the norm-distance from $A$ to $\cB$. One motivation 
for studying such quotient seminorms is that they arise naturally when considering
matricial Lipschitz seminorms, as defined in \cite{Wuw1, Wuw2, Wuw3} (where
the Leibniz aspect was not used)
and as discussed in the
final section of \cite{R21}. These involve a unital C*-algebra $\cD$, and for
each natural number $m$ the corresponding matrix C*-algebra 
$\cA = M_m(\cD) = M_m(\bC) \otimes \cD$ over $\cD$. The 
Leibniz seminorm on $\cA$ from a matricial
Leibniz seminorm on $\cD$ will take value 0 on the C*-subalgebra
$\cB = M_m(\bC) \otimes 1_\cD$, and so can be viewed as a seminorm
on $\cA/\cB$. 
Furthermore, in this situation $\cB$ is finite-dimensional,
and so elements of $\cA$ will always have a best approximation by an element
of $\cB$ (probably not unique). This is good motivation for studying best 
approximations. Although there is a 
very large literature dealing with best approximation of elements of a Banach
space by elements of a closed subspace \cite{Sng}, I have found almost no literature
concerning this topic for the case of C*-algebras and their C*-subalgebras. 
The later sections of this paper are devoted to developing some basic results
about this for C*-algebras. We show that the situation for C*-algebras has
some nice features in comparison to the general case, and we are able
to tie together a few results scattered in the literature.

Somewhat to my surprise I have found that the seminorms from quotients
of C*-algebras are usually strongly-Leibniz. To be specific, the most important
theorem of this paper states:

\begin{thmnonum}
\label{thmMainintro} 
Let $\cA$ be a C*-algebra and let $\cB$ be a C*-subalgebra of $\cA$. 
Assume that $\cB$ contains a bounded approximate identity for $\cA$.
Let $L$ be defined as above. Then $L$ is Leibniz, that
is,
\[
L(AC) \leq L(A)\|C\| + \|A\|L(C)
\]
for all $A, C \in \cA$. If $\cA$ is unital and if $1_\cA \in \cB$, then $L$ is
strongly-Leibniz, that is, $L$ is Leibniz, $L(1) = 0$, and
if $A$ is invertible in $\cA$ then
\[
L(A^{-1}) \leq \|A^{-1}\|^2L(A)  .
\]
\end{thmnonum}

This theorem is a fairly simple consequence of a C*-algebra version of the
Arveson distance formula \cite{Arv} that has been widely used in the study of nest algebras. This C*-algebra version and the proof of the above theorem are given in Section \ref{secC^*}.

One can, of course, ask what happens for Banach algebras $\cA$ and their closed
subalgebras $\cB$. In Section \ref{sec1} we discuss a property, that we call the
``same-norm approximation property'', that guarantees that the seminorm on
$\cA$ pulled back from the quotient norm on $\cA/\cB$ is Leibniz. We then show
that if $\cA$ is a C*-algebra and $\cB$ is a C*-subalgebra of $\cA$ that is
\emph{central} in $\cA$, then $\cB$ has the same-norm approximation property
in $\cA$.

But a main application of our later discussion of best approximation is to
give an example of a finite-dimensional C*-algebra and a (non-central)
unital C*-subalgebra that does not have the same-norm approximation property. 
The nature
of the construction strongly suggests that the same-norm approximation
property will usually fail unless one is dealing with central subalgebras.

It is a pleasure to thank Man-Duen Choi, Erik Christensen, and Dick Kadison
for very helpful conversations concerning the subject of this paper.


\section{The Leibniz inequality and Banach algebras}
\label{sec1}

In this section we make some elementary observations in the setting of Banach algebras.
These observations relate the Leibniz inequality for quotient
seminorms with a property of approximations by elements of a subalgebra. We then
apply these observations to the setting of C*-algebras. 

Our experience with Hilbert spaces leads us to expect that if $B$ is a best
approximation to $A$ then $\|B\| \leq \|A\|$. But this easily fails for the seminorms coming
from best approximation in C*-algebras, as we will see later. In general, the most that one
can say is that 
\[
\|B\| \leq 2\|A\|
\]
because $B$ must be at least as close to $A$ as is $0$. Easy examples show that 2 is the
best constant that holds here in general. But the stronger inequality is desirable for our purposes,
as we will show in this section, and so we make the following definition:

\begin{definition}
\label{samenorm}
Let $\cA$ be a normed vector space and let $\cB$ be a closed subspace
of $\cA$. Let $M$ be the quotient norm on $\cA/ \cB$, pulled back to $\cA$.
That is,
\[
M(A) = \inf \{\|A-B\|: B \in \cB\}  .
\]
Given $A \in \cA$, we say that $A$ is \emph{same-norm approximatable}
in $\cB$ if for every $\e > 0$ there is a $B \in \cB$ with $\|B\| \leq \|A\|$ and
$\|A-B\| < M(A) + \e$. If every element of $\cA$ is same-norm 
approximatable in $\cB$
then we say that $\cB$ has the \emph{same-norm approximation property} in $\cA$.
\end{definition}

\begin{proposition}
 \label{same-leib}
 Let $\cA$ be a normed algebra, let $\cB$ be a closed subalgebra
 of $\cA$, and let $M$ be the quotient norm on $\cA/ \cB$, pulled back to $\cA$.
 Let $A \in \cA$. If $A$ is same-norm approximatable in $\cB$, then for every
 $C \in \cA$ the Leibniz inequality
 \[
 M(AC) \leq M(A)\|C\| +\|A\|M(C)
\]
holds for all $C \in \cA$, and similarly for $M(CA)$. 
  \end{proposition}
  
\begin{proof}
Given $\e > 0$, choose $D \in \cB$ such that $\|C-D\| < M(C) + \e$,
and choose $B\in \cB$ such that $\|B\| \leq \|A\|$ and $\|A-B\| < M(A) + \e$.
Then
\begin{align*}
M(AC) &\leq \|AC - BD\| \leq \|A-B\|\|C\| + \|B\|\|C-D\|   \\
& \leq (M(A) + \e)\|C\| + \|A\|(M(C) + \e).
\end{align*}
Since $\e$ is arbitrary, we obtain the desired inequality.
A similar calculation works for $M(CA)$.
\end{proof}  

We remark that the above proof is somewhat parallel
to the proof of proposition 5.2 of \cite{R21}.

In the setting of normed algebras it would be interesting to 
have conditions that also ensure that $M$ is strongly-Leibniz,
as defined earlier.
For this to make sense we need that $\cA$ is unital, (There
may well be a suitable generalization of the notion of
strongly-Leibniz to non-unital algebras, but I have not
explored this possibility.)

For C*-algebras the main situation in which the same-norm
approximation property seems to arise is the following.

\begin{proposition}
\label{prop-comm}
Let $\cA$ be a C*-algebra and let $\cB$ be a C*-subalgebra
of $\cA$ such that $\cB$ is \emph{central} in $\cA$. Then $\cB$ has
the same-norm approximation property in $\cA$.

More specifically, let $A \in \cA$, let $F \in \cB$, and let $G$ be the
radial retraction of $F$ to the ball about 0 of radius $\|A\|$ in $\cB$
defined by the continuous-function calculus for normal elements. Then 
$\|A-G\| \leq \|A-F\|$, and, of course, $\|G\| \leq \|A\|$.
\end{proposition}

\begin{proof}
The assertion in the second paragraph implies the assertion in the first
paragraph because it shows that whatever approximation we have,
we can always get from it another that is as close but also
satisfies the same-norm condition. 

If $\cA$ is not unital, then when we adjoin an identity element to
$\cA$ the subalgebra $\cB$ is still central, and distances are
not changed. So we will assume now that
$\cA$ is unital. Let $\tilde \cB$ denote $\cB$ with the identity element 
of $\cA$ adjoined (if it is not already in $\cB$).

To prove the assertion of the
second paragraph, notice that 
when we view $\tilde \cB$ as the algebra of continuous functions
on a compact space, we see that $G = HF$ where $H = K^{-1}$
and $K = \max\{1, |F|/ \|A\|\}$ as functions. Notice that $K - 1$ will be in
$\cB$, as will then $H-1$, so that $K$ and $H$ are both multipliers
of $\cB$. To prove the assertion of the second paragraph,
it suffices to show that as operators
\[
(A-F)^*(A-F) \ \geq \ (A-G)^*(A-G)  ,
\]
that is, upon simplification, that 
\[
|F|^2(1-H^2) \ \geq \ (1-H)(F^*A + A^*F)  .
\]
Let $(\pi, \cH)$ be an irreducible representation of $\cA$, 
so that $\pi(F)$ is a scalar multiple, say, $\l$, of $I_\cH$, since $\cB$ is central. 
Then $\pi(H) = \mu I_\cH$ where $\mu = 1$ if $|\l| \leq \|A\|$ while 
$\mu = \|A\|/|\l|$ otherwise. Then for the above inequality we need to know that
\[
|\l|^2(1-|\mu|^2) \geq (1-\mu)(\bar \l \pi(A) + \l \pi(A^*)).
\]
Now $\bar \l \pi(A) + \l \pi(A^*)$ has norm $\leq 2|\l|\|A\|$ and is Hermitian,
and thus it suffices to know that
\[
|\l|^2(1-|\mu|^2) \geq (1-\mu)2|\l|\|A\|   .
\]
If $\mu = 1$ then both sides are 0, while if $|\l| > \|A\|$ so that 
$\mu = \|A\|/\l < 1$ then a simple calculation shows that this inequality
holds. Since all of this works for any irreducible representation, we
obtain the desired inequality.
\end{proof}

\begin{corollary}
\label{central}
Let $\cA$ be a C*-algebra and let $\cB$ be a C*-subalgebra
of $\cA$ such that $\cB$ is \emph{central} in $\cA$, and let 
$M$ be the quotient norm on $\cA/\cB$ pulled back to $A$. 
Then $M$ satisfies the Leibniz inequality.
\end{corollary}

This corollary is related to the main result of \cite{Gjn}, which, however,
is more in the spirit of the next section.

We remark that if $\cA$ is unital
and if $1_\cA \notin \cB$, then $M(1_\cA) \neq 0$,
so that $M$ can not be strongly-Leibniz.

\begin{corollary}
\label{identity}
Let $\cA$ be a unital C*-algebra, let $\cB = \bC 1_\cA$, and let $M$ be
the corresponding seminorm on $\cA$. Then $M$ satisfies the Leibniz
inequality.
\end{corollary}

This last corollary shows that we do obtain a Leibniz seminorm 
in the $m = 1$  case of the situation
in which $\cA = M_m(\cD)$ and $\cB = M_m(\bC 1_\cD)$
that was mentioned in introduction. But it does not show that
the seminorm is strongly-Leibniz. (We will see in the next section
that it is indeed strongly-Leibniz.)
This corollary is also strongly related to an often-cited result of Stampfli \cite{Stm},
and to, for example, equation 1.6 of \cite{BhS} and section 3.1
of \cite{ChL}.

We remark that if $\cA = M_2(\bC)$ and if $\cB$ is the C*-subalgebra of
all diagonal matrices in $\cA$, then $\cB$ has the same-norm approximation
property in $\cA$. This is seen easily by direct calculation. Thus a proper C*-subalgebra
does not need to be central in order to have the same-norm approximation property.
But it would be interesting to know if there are any other C*-algebras that
have a proper non-central C*-subalgebra that has the same-norm approximation
property. This may be related to the fact shown in \cite{Pd2} that $M_2(\bC)$
is the only C*-algebra that has a proper Chebyshev C*-subalgebra of dimension
strictly greater than 1.

\begin{corollary}
\label{normal}
Let $\cA$ be a C*-algebra, let $\cB$ be a C*-subalgebra
of $\cA$, and let $M$ be the quotient norm on $\cA/\cB$ pulled back to $A$.
Let $A\in \cA$. If there is a best approximation $B \in \cB$ to $A$ which is normal 
and commutes with $A$, then $A$ is same-norm approximatable by $\cB$
(so Proposition \ref{same-leib} is applicable).
\end{corollary}

\begin{proof}
By Fuglede's theorem (theorem 4.76 of \cite{Dgl}) $A$ also commutes
with $B^*$. Thus the C*-algebra $C^*(B)$ generated by $B$
is a central C*-subalgebra of the
C*-algebra $C^*(A, B)$ generated by
$A$ and $B$. Thus we can apply Proposition \ref{prop-comm}.
Since $C^*(B)$ contains a best approximation  to
$A$ in $\cB$, namely $B$, the desired conclusion follows.
\end{proof}

\begin{corollary}
\label{commutator}
Let $\cD$ be a C*-algebra, and for a natural number $m$ 
let $\cA = M_m(\cD)$ and let $\cB$ be the 
C*-subalgebra $\cB = M_m(\bC) \otimes 1_\cD$
as mentioned in the introduction.  Then any element $A$ of $\cA$ that is
in the commutant of $\cB$ in $\cA$ (i.e. is in $I_m \otimes \cD$)
is same-norm approximatable in $\cB$ (so that Proposition \ref{same-leib}
is applicable).
\end{corollary}

\begin{proof}
It is easily seen that if $A = I_m\otimes D$ for some $D \in \cD$, and if
$\l \in \bC$ is such that $\l 1_\cD$ is a best approximation to $D$ in $\bC 1_\cD$, chosen via Corollary \ref{identity} so that $|\l | \leq \|D\|$, 
then $B = \l I_m \otimes 1_\cD$ is a best approximation to $A$ in $\cB$
and $\|B\| \leq \|A\|$.
\end{proof}

It would be interesting to have examples of Banach algebras that
are not C*-algebras but that have subalgebras
that satisfy the same-norm approximation property.


\section{Many quotient norms of C*-algebras  are \\
strongly Leibniz}
\label{secC^*}

Let $\cA$ be a C*-algebra and let $\cB$ be a C*-subalgebra of $\cA$. 
As before we
let $L$ denote the pull-back to $\cA$ of the corresponding quotient norm
on $\cA/\cB$, so that 
\[
L(A) = \inf \{\|A-B\|: B \in \cB\}  
\]
for all $A \in \cA$. 
For the purposes of my general investigation of Leibniz
seminorms on C*-algebras, the most important theorem of this paper is the following:

\begin{theorem}
\label{thmMain}
Let $\cA$ be a C*-algebra and let $\cB$ be a C*-subalgebra of $\cA$. 
Assume that $\cB$ contains a bounded approximate identity for $\cA$.
Let $L$ be defined as above. Then $L$ is Leibniz, that
is,
\[
L(AC) \leq L(A)\|C\| + \|A\|L(C)
\]
for all $A, C \in \cA$. If $\cA$ is unital and if $1_\cA \in \cB$, then $L$ is
strongly-Leibniz, that is, $L$ is Leibniz, $L(1) = 0$, and
if $A$ is invertible in $\cA$ then
\[
L(A^{-1}) \leq \|A^{-1}\|^2L(A)  .
\]
\end{theorem}

In proposition 1.2iii of \cite{R21} it is seen that any supremum of
a (possibly infinite) family of Leibniz seminorms will again be
Leibniz, with a similar statement for strongly-Leibniz seminorms. 
In proposition 2.1 of \cite{R21} it is seen
that for any derivation $d$ from $\cA$ into a normed bimodule, if
we set $L_d(A) = \|d(A)\|$ then $L_d$ is a Leibniz seminorm
on $\cA$, and if $\cA$ is unital and the bimodule is non-degenerate,
then $L$ is strongly-Leibniz. Combining these two facts gives the last part of
the next theorem. Furthermore, Theorem \ref{thmMain} is 
an immediate consequence
of this next theorem.

\begin{theorem}
\label{thmArv}
Let $\cA$ be a C*-algebra and let $\cB$ be a C*-subalgebra of $\cA$.
Assume that $\cB$ contains a bounded approximate identity for $\cA$.
Let $L$ be defined as above.
For any $A \in \cA$ there is a non-degenerate $*$-representation,
$(\cH, \pi)$, of $\cA$, and a Hermitian unitary operator $U \in \cL(\cH)$,
such that $[U, \pi(B)] = 0$ for all $B \in \cB$, and
\[
L(A) = (1/2)\|[U, \pi(A)]\|   .
\]
Furthermore, $L$ is the supremum of all the seminorms $L_{(\cH, \pi, U)}$
defined by $L_{(\cH, \pi, U)}(C) = (1/2)\|[U, \pi(C)]\|$ for $C \in \cA$,
as $(\cH, \pi, U)$ ranges over all 
non-degenerate $*$-representations, $(\cH, \pi)$,
of $\cA$ and all Hermitian unitary elements $U \in \cL(\cH)$ such that
$[U, \pi(B)] = 0$ for all $B \in \cB$. Since each $L_{(\cH, \pi, U)}$
is Leibniz, so is $L$.  If $\cA$ is unital and if $1_\cA \in \cB$, then $L$ is
strongly-Leibniz.
\end{theorem}

Theorem \ref{thmArv} is in turn a C*-algebraic 
variation on Arveson's distance
formula for nest algebras \cite{Arv}, 
as used in the case of von Neumann algebras
\cite{Chr}. I thank Erik Christensen for elucidating for me this 
formula for the case of 
von Neumann algebras.

I have not seen a good way to weaken the assumption
that $\cB$ contains a bounded approximate identity for $\cA$.

\begin{proof}[Proof of theorem \ref{thmArv}]
Notice that if $C\in \cA$ then for any $(\cH, \pi, U)$
as above we have
$\|[U, \pi(C)]\|=\|[U, \pi(C-B)]\| \leq 2\|C-B\|$
for all $B\in \cB$, and thus 
\[
L_{(\cH, \pi, U)}(C) \leq L(C)
\]
for all $C \in \cA$. The fact that $L$ is the supremum
of the seminorms $L_{(\cH, \pi, U)}$ thus follows from
the first statement in the theorem. The fact that $L$ is
strongly-Leibniz if $\cA$ is unital and if $1_\cA \in \cB$
then follows from the comments made
before the statement of the theorem. Thus it remains to
prove the first statement of the theorem.

So let $A \in \cA$ be given. By scaling, we see that it suffices
to treat the case in which $L(A) = 1$, and
so for ease of bookkeeping we assume this. We now use
the first basic tool of linear approximation theory. By applying
the Hahn-Banach theorem to the image of $A$ in $\cA/\cB$,
we see that there is a linear functional, $\psi$, on $\cA$,
such that $\|\psi\| = 1$, $\psi(A) = 1$, and $\psi(B)=0$
for all $B \in \cB$. Our proof is then based on the following
key lemma, which is very closely related to the polar decomposition
of linear functionals in the predual of a von Neumann algebra. See
notably corollary 7.3.3 of \cite{KR2}. This lemma can also be 
obtained by examining the universal representation of $\cA$,
as suggested in the proof of lemma 1 of \cite{Pd2}.

\begin{lemma}
\label{lemKey}
Let $\cA$ be a C*-algebra, and let $\psi$ be a linear functional
on $\cA$ such that $\|\psi\| = 1$. Then there is a cyclic
$*$-representation, $(\cH, \pi, \xi)$, of $\cA$, and a vector $\eta \in \cH$,
such that $\|\xi\| = 1 = \|\eta\|$ and
\[
\psi(C) = \<\pi(C)\xi, \eta\>
\]
for all $C \in \cA$.
\end{lemma}

Before giving the proof of this lemma, we show how to use it to 
complete the proof of Theorem \ref{thmArv}. As in the first fragment
of the proof of Theorem \ref{thmArv} given above, let $\psi$
be the linear functional on $\cA$ with $\|\psi\| = 1$, $\psi(A) = L(A) = 1$,
and $\psi(B) = 0$ for all $B \in \cB$. According to Lemma \ref{lemKey}
there is a non-degenerate $*$-representation  $(\cH, \pi)$, of $\cA$, 
and vectors $\xi, \eta \in \cH$,
such that $\|\xi\| = 1 = \|\eta\|$ and
\[
\psi(C) = \<\pi(C)\xi, \eta\>
\]
for all $C \in \cA$. Let $P$ denote the orthogonal projection
of $\cH$ onto the closure of $\pi(\cB)\xi$. Since the closure of 
 $\pi(\cB)\xi$ is clearly $\cB$ invariant, we have $[P, \pi(B)] = 0$
 for all $B \in \cB$. Notice that for all $B \in \cB$ we have
$0 = \<\pi(B)\xi, \eta\>$, and thus $\eta \perp \pi(\cB)\xi$, so that
 $P(\eta) = 0$. Also, because $\cB$ contains a bounded 
 approximate identity for $\cA$, we see that $\xi$
 is in $\pi(\cB)\xi$, so that $P(\xi) = \xi$. 
 For notational simplicity, let us now write $C$ for $\pi(C)$, etc. Then
 \[
 \<[C, P]\xi, \eta\> = \<CP\xi, \eta\> - \<PC\xi, \eta\> = \<C\xi, \eta\> = \psi(C)  
 \]
 for all $C \in \cA$. In particular, $\|[A, P]\| \geq |\psi(A)| = L(A)$. 
 This does not quite fit our needs, as
 we only have the general estimate that for $C \in \cA$ and $B \in \cB$ 
 \[
 \|[C, P]\| = \|[C-B, P]\| \leq 2\|C-B\|   ,
 \]
so that $\|[C, P]\| \leq 2L(C)$, and we do not want the factor of 2 here.
(But notice the importance of having P commute with all the $B$'s.)
To rectify this, let $U = 2P - I$, so that $U$ is a Hermitian unitary, and
in particular, $\|U\| = 1$. Then we will have
\[
\psi(C) = (1/2)\<[C, U]\xi, \eta\>   
\]
for all $C \in \cA$, 
so that in particular $(1/2)\|[U, A]\| \geq L(A)$. But for any $C \in \cA$ and $B \in \cB$
we still have, much as above,
\[
\|[C, U]\| = \|[C-B, U]\| \leq 2\|C-B\|,
\]
so that $(1/2)\|[U, C]\| \leq L(C)$. This gives the desired result.
\end{proof}

For a related result in the finite-dimensional case see remark 3.1
of \cite{BhS}.

\begin{proof}[Proof of Lemma \ref{lemKey}]
For the reader's convenience we now give a proof 
of Lemma \ref{lemKey} that uses as its biggest tool just the Jordan 
decomposition of Hermitian linear functionals into differences of 
two positive linear functionals, 

Let $\cA$ and $\psi$ be as in the statement of Lemma \ref{lemKey}.
If $\cA$ is not unital, adjoin an identity element in the usual way so
as to obtain a C*-algebra. By the technical step in the proof of the
Hahn-Banach theorem we can extend $\psi$ to the resulting C*-algebra,
with no increase in its norm. Until the end of the proof we will now assume
that $\cA$ is unital.

Form the algebra $M_2(\cA)$ of $2 \times 2$ matrices with entries in $\cA$,
with its unique C*-algebra structure. Define on $M_2(\cA)$ a linear functional,
$\psi_2$, by
\[
\psi_2\left( \begin{pmatrix}  A & C \\
                                          B & D       \end{pmatrix}    
                  \right) 
= (1/2)(\psi(B) \ + \overline{ \psi (C^*)})              .         
\]
Then $\psi_2$ is a Hermitian linear functional, that is, 
$\psi_2(T^*) = \overline{\psi(T)}$
for any $T \in M_2(\cA)$. Note also that $\|\psi_2\| = 1$.

Let $\psi_2 = (1/2)(\phi^+ - \phi^-)$ be the Jordan decomposition
of $\psi_2$
(theorem 4.3.6 and remark 4.3.12 of \cite{KR1}), so that $\phi^+$
and $\phi^-$ are 
positive linear functionals on $M_2(\cA)$ such that 
$\|\psi^+\| \ + \ \|\phi^-\| = 2 \|\psi_2\| = 2$. Note that $\psi_2(I_2) = 0$,
where $I_2$ is the identity element of $M_2(\cA)$. It follows
that $\phi^+(I_2) = \phi^-(I_2)$, so that $\|\phi^+\| = \|\phi^-\| = 1$.
Thus both $\phi^+$ and $\phi^-$ are states on $M_2(\cA)$.

Let $(\cH^+, \pi^+, \xi^+)$ and $(\cH^-, \pi^-, \xi^-)$ be the GNS 
representations for $\phi^+$ and $\phi^-$. Set
\[
\cH_2 = \cH^+ \oplus \cH^-  \quad \quad \mathrm{and} \quad \quad \pi_2 
= \pi^+ \oplus \pi^-  ,
\]
and set
\[
\xi_2 = (\xi^+ \oplus \xi^-)/\sqrt 2 \quad \quad \mathrm{and} \quad \quad
\eta_2 = (\xi^+ \oplus -\xi^-)/\sqrt 2.
\]
Note that $\|\xi_2\| = 1 = \|\eta_2\|$. Then for any $T \in M_2(\cA)$ we
have
\begin{align*} 
\psi_2(T) &= (1/2)(\phi^+(T) - \phi^-(T))  \\
&= (1/2)(\<\pi^+(T)\xi^+, \xi^+\> - \<\pi^-(T)\xi^-, \xi^-\>)  \\
&= \<\pi_2(T)\xi_2, \eta_2\>  .
\end{align*} 

Now for any $A \in \cA$ we have
\[
\begin{pmatrix}  0 & 0 \\
                          A & 0       \end{pmatrix}    
=
\begin{pmatrix}  0 & 1 \\
                          1 & 0       \end{pmatrix}    
                                          \begin{pmatrix}  A & 0 \\
                                          0 & 0       \end{pmatrix}    .                                                    
\]
Thus
\[
\psi(A) = \psi_2\left( \begin{pmatrix}  0 & 0 \\
                                                          A & 0       \end{pmatrix} \right) 
= \<\pi_2 \left( \begin{pmatrix}  A & 0 \\
                                          0 & 0       \end{pmatrix}    \right) \xi_2, \ \
\pi_2 \left( \begin{pmatrix}  0 & 1 \\
                                          1 & 0       \end{pmatrix}    \right) \eta_2 \>  . 
\]
Let $P = \pi_2(
(\begin{smallmatrix}  1 & 0 \\ 0 & 0 \end{smallmatrix} ) ) $,
and then set $\cH = P\cH_2$. Notice that $\pi_2(
(\begin{smallmatrix}  A & 0 \\ 0 & 0 \end{smallmatrix} ) ) $ carries
$\cH$ into itself for every $A \in \cA$. We obtain in this way a unital
*-representation, $\pi$, of $\cA$ on $\cH$. Set $\xi = P\xi_2$ and $\eta =
P\pi_2( (\begin{smallmatrix}  0 & 1 \\ 1 & 0 \end{smallmatrix} ) )\eta_2 . $
From the equation displayed just above we see that
\[
\psi(A) = \<\pi(A)\xi, \eta\>
\]
for all $A \in \cA$. From the definitions given above it is easily
seen that $\|\xi\| \leq 1$ and $\|\eta\| \leq 1$. 
One can restrict $\pi$ to the cyclic
$\cA$-subspace of $\cH$ generated by $\xi$. Then from the fact
that $\|\psi\| = 1$ it follows that $\|\xi\| = 1$ and $\|\eta\| = 1$,
as desired. If an identity element had been adjoined to the 
original $\cA$, one can restrict $\pi$ further to the cyclic
$\cA$-subspace of $\cH$ generated by $\xi$.
\end{proof}

If we now let $\{\cH, \pi\}$ be the Hilbert-space direct sum over all
$A \in \cA$ of the representations obtained in Theorem \ref{thmArv},
and if we let $U$ be the direct sum of the corresponding Hermitian
unitaries obtained there, than we immediately obtain:

\begin{corollary}
\label{corsingle}
Let $\cA$ be a C*-algebra, and let $\cB$ be a C*-subalgebra of $\cA$ that  
contains a bounded approximate identity for $\cA$.
Let $L$ be defined as above for $\cA$ and $\cB$. Then there is a non-degenerate unitary
representation $\{\cH, \pi\}$ of $\cA$ and a Hermitian unitary operator
$U$ on $\cH$ that commutes with $\pi(B)$ for every $B \in \cB$, such that
\[
L(A) = (1/2)\| [U, A]\|
\]
for every $A\in \cA$.
\end{corollary}

Note that if $L$ is defined as above, and if $M$ is any seminorm on $\cA$ that satisfies the conditions
that $M(A) \leq \|A\|$ for all $A \in \cA$, and 
that $M(B) = 0$ for all $B \in \cB$,
then we have $M(A) \leq L(A)$ for all $A \in \cA$. Thus we obtain 
immediately the following corollary, which is of interest in particular
for the situation discussed in the introduction
in which we have a unital C*-algebra $\cD$ and 
$\cA = M_n(\bC) \otimes \cD$ while $\cB = M_n(\bC) \otimes 1_\cD$.

\begin{corollary}
\label{cormax}
Let $\cA$ be a C*-algebra, and let $\cB$ be a C*-subalgebra of $\cA$
that contains a bounded approximate identity for $\cA$.
Let $L$ be defined as above. Then $L$ is the maximal Leibniz
seminorm on $\cA$ that takes value 0 on $\cB$ and is dominated by
the C*-norm of $\cA$. If $\cA$ is unital and if $1_\cA \in \cB$, then $L$ is
the maximal strongly-Leibniz seminorm on $\cA$ that takes 
value 0 on $\cB$ and is dominated by
the C*-norm of $\cA$.

\end{corollary}

Presumably Arveson's distance formula provides strongly-Leibniz
seminorms on suitable nest algebras, but I have not explored this matter.

We remark that from the proof given above for Theorem \ref{thmArv} one can
extract a corollary that is closely related to lemma 2.1 of \cite{HLT} (which deals
with finite-dimensional matrix algebras).

Now let the notation be as in Corollary \ref{corsingle}, and let $V = iU$,
so that $V$ is a skew-Hermitian operator. Define $\d_V$ on $\cA$
by $\d_V(A) = (1/2)[V, A]$. Then $\d_V$ is a $*$-derivation
of $\cA$ on $\cH$, where the general definition is:

\begin{definition}
Let $\cA$ be a C*-algebra, and let $\cH$ be a Hilbert space. 
We say that a pair $(\pi, \d)$ is a \emph{$*$-derivation of $\cA$ on $\cH$}
if $\pi$ is a non-degenerate $*$-representation of $\cA$ on $\cH$
and $\d$ is a linear function from $\cA$ into $\cB(\cH)$ such that
$\d(A^*) = (\d(A))^*$ and 
\[
\d(AC) = \d(A)C + A\d(C)
\]
for all $A, C \in \cA$.
By the \emph{kernel of $(\pi, \d)$} we will mean the kernel of $\d$.
\end{definition}

It is known that every $*$-derivation on a C*-algebra is necessarily
continuous. (See 4.6.66 of \cite{KR1} and  4.6.66 of \cite{KR3}.)
It is evident that the kernel of a $*$-derivation 
is a C*-subalgebra of $\cA$. If $\cA$ is
unital then it is easily seen that $\d(1) = 0$, so that the kernel
of $(\pi, \d)$ is a unital C*-subalgebra of $\cA$. Then Corollary
\ref{corsingle} quickly gives, by means of $\d_V$, the following characterization of the kernels
of $*$-derivations of unital C*-algebras:

\begin{corollary}
Let $\cA$ be a unital C*-algebra. Then the possible kernels of
$*$-derivations of $\cA$ on various Hilbert spaces are exactly the
unital C*-subalgebras of $\cA$.
\end{corollary}  

If $\cA$ is not unital, then I do not know whether the kernel of
a $*$-derivation must contain a bounded approximate identity
for $\cA$, and thus I do not know how to characterize all of the
possible kernels of $*$-derivations of $\cA$.


\section{Chebyshev subalgebras}
\label{sec2}

We now begin our study of best approximations in C*-algebras. We start  by recalling some relevant well-known facts about
approximation of Banach-space elements by elements of
a closed subspace, and by relating them to our C*-algebra situation.
Much of the literature on that subject seems to be motivated by
Chebyshev's famous theorem that a best uniform approximation of a continuous
real-valued function on the unit interval by polynomials of degree no 
greater than a given natural number is unique. (See, e.g., section 7.6 of \cite{Dvs}.)
Accordingly, when
$\cA$ is a Banach space and $\cB$ is a closed subspace, one says
that $\cB$ is a \emph{Chebyshev subspace} if every element of $\cA$
has a best approximation by an element of $\cB$, and that best
approximation is unique. For example, every closed subspace
of a Hilbert space is a Chebyshev subspace.

Following partial results by A. G. Robertson \cite{Rbt, RbY} and a few others, Gert
Pedersen showed \cite{Pd2} that if $\cA$ is a unital C*-algebra
and $\cB$ is a unital C*-subalgebra that is a Chebyshev subspace
of $\cA$, then either $\cB = \cA$, or $\cB = \bC 1_\cA$, or 
$\cA = M_2(\bC)$ and $\cB$ is the subalgebra of diagonal
matrices.

\begin{example}
\label{non-unique}
Let $\cA$ be $(M_2(\bC))^3$, the C*-algebra of 3-tuples
of elements of $M_2(\bC)$, 
and let $\cB$ be its C*-subalgebra of constant 3-tuples. Let
\[
A = \left\{ \begin{pmatrix}  1 & 0 \\
                                          0 & 1       \end{pmatrix}     \ , \ 
               \begin{pmatrix}  1 & 0  \\
                                         0 & -1             \end{pmatrix}     \ , \
              \begin{pmatrix}   1 & 0  \\                                              
                                         0 & 0             \end{pmatrix}     
     \right\}  .         
\]     
Since the distance between the first two entries of $A$ is 2, it is
clear that any $B =   \begin{pmatrix}   t & 0  \\                                              
                                                                 0 & 0             \end{pmatrix}    
$   
(viewed as a constant 3-tuple) with $0 \leq t \leq 2$ is a best   
approximation to $A$ in $\cB$. In particular, for $t=2$ we see that
$B$ is a best approximation with $\|B\| > \|A\|$. But since $B$ for
$t=0$ is also a best approximation, $A$ is still same-norm
approximatable in $\cB$. In the last section of this
paper we will give examples
that are not same-norm approximatable but have a unique
best approximation.
\end{example}

We remark that anyway, in the Banach space case, same-norm
approximation can easily fail for Chebyshev subspaces. This
already happens for uniform approximation of continuous
real-valued functions on the unit interval by polynomials, as is 
seen when uniformly approximating the function $f(t) = t^{1/n}$ by polynomials
of degree 1 or less. As $n$ goes to $\infty$ the norm of the best approximating
linear polynomial goes to 2. 

However, one striking phenomenon concerning polynomial approximation 
of continuous real-valued functions on the interval is that the error function 
equi-oscillates. (See, e.g., theorem 7.6.2 of \cite{Dvs}.) 
That is, if $f \in C_\bR([0,1])$ and if $p$
is the best uniform approximation to $f$ by polynomials of degree no greater than
$n$, then one can find $n+2$ distinct points in $[0,1]$ at which $|f-p|$
takes its maximum value with the sign of $f-p$ evaluated at these
points alternating at successive points. We will see an echo of this 
equi-oscillation later in our C*-algebraic setting. 
                                  

\section{Witnesses for best approximation in Banach spaces}
\label{wit-Ban}

If we want to give specific counter-examples involving best
approximations, we need to prove that a candidate best approximation is
indeed a best approximation. For this purpose 
we now discuss further the first basic tool of linear approximation theory
\cite{Sng} that we used already in the second paragraph of the
proof of Theorem \ref{thmArv}. In the next sections 
we will see that in our C*-algebraic 
setting this tool has some nice special features. 

Let $\cA$ be a Banach space, let $\cB$ be a proper 
closed subspace of $\cA$, and let $M$
be the quotient norm on $\cA / \cB$, often viewed as a seminorm on $\cA$.
Let $A \in \cA$. Denote by $\cA'$ the Banach-space dual of $\cA$, and
by $\cAo$ its unit ball. Then, much as in the proof of Theorem \ref{thmArv}, 
by the Hahn-Banach theorem applied to
the image of $A$ in $\cA / \cB$ there is a $\psi \in \cAo$ such that $\|\psi\| = 1$,
$\psi \in \cB^\perp$ (in the sense that $\psi(B) = 0$ for all $B \in \cB$)
and $\psi(A) = M(A)$. It may well happen that $\psi$ is not unique.

Suppose that $A$ has a best approximation, $B$, in $\cB$. Then
$\|A-B\| = M(A) = \psi(A)$. The following proposition is a basic
tool for proving that one has a best approximation.

\begin{proposition}
\label{prop-wit}
Let $\cA$, $\cB$, and $M$ be as above, and let $A \in \cA$ and
$B \in \cB$ be given. Let $\psi \in \cAo$ be given. If $\psi \in \cB^\perp$
and $\psi(A) = \|A-B\|$, then $B$ is a best approximation to $A$ in $\cB$.
\end{proposition}

\begin{proof}
Note that necessarily $\|\psi\| = 1$, as long as $A \notin \cB$. For any $C \in \cB$ we have
\[
\psi(A - C) = \psi(A) = \|A-B\|  .
\]
It follows that $\|A-C\| \geq \|A-B\|$.
\end{proof}

Thus we can say that $\psi$ is a ``witness'' for the fact that $B$
is a best approximation to $A$ in $\cB$. That is, if we can find such
a $\psi$ then we have a proof that $B$ is a best approximation.

As suggested by our comments above concerning equi-oscillation,
it is useful to look at the error of an approximation, that is, at $Z = A-B$.
If $B$ is a best approximation to $A$ in $\cB$, then it is clear that $0$
is a best approximation to $Z$ in $\cB$ (thus our use of the symbol
$Z$, but see also its use in \cite{DMR, AMM}), and if $\psi$ is a witness for $B$ as
above, then we have $\psi(Z) = \psi(A) = \|Z\|$. 

\begin{definition}
\label{def-min}
With $\cA$ and $\cB$ as above, we say that an element $Z \in \cA$
is \emph{minimal for $\cB$}, or \emph{$\cB$-minimal}, if $0$ is a best approximation to $Z$ 
in $\cB$, so that $M(Z)=\|Z\|$, and $M(Z+B) = \|Z\|$ for every
$B \in \cB$.
\end{definition} 

Then a witness to the $\cB$-minimality of $Z$ will be a $\psi \in \cAo$
such that $\psi \in \cB^\perp$ and $\psi(Z) = \|Z\|$. Note that for
any $B \in \cB$ we will have 
\[
M(Z+B) = \|Z\| = \psi(Z).
\]
This can be of help in calculating $M(A)$ for $A \in \cA$.

Suppose now that $\cB$ is finite-dimensional, so that best approximations
always exist. If $\cA$ is over the complex numbers, we will forget that and
view $\cA$ as being over $\bR$. This has no effect on best approximations, 
but dimensions will usually be over $\bR$ as this simplifies a bit the bookkeeping.

Recall Caratheodory's theorem 
(see, e.g., exercise 19 of chapter 3 of \cite{Rdn}) that if $\cC$ 
is a finite-dimensional Banach space of dimension $q$
and if $K$ is a closed bounded convex subset of $\cC$, with $K_e$
its set of extreme points, then $K=\mathrm{convex}(K_e)$, with no need
to take closure, and in fact that every element of $K$ is a convex
combination of at most $q+1$ elements of $K_e$. Furthermore,
one can show that if the element is in the
boundary of $K$ then at most $q$ elements of $K_e$ are needed.
(See lemma 1.1 of chapter II of \cite{Sng}.)

Suppose now that $A \in \cA$ but $A \notin \cB$. Let $\cF = \cB \oplus \bR A$,
a subspace of $\cA$, equipped with the norm from $\cA$. Let 
$p = \dim_\bR(\cB)$, so that $\dim_\bR(\cF) = p+1$. When $A$ is viewed as
an element of $\cF$, its best approximations in $\cB$ are the same as
the best approximations when $A$ is viewed as an element of $\cA$.
Suppose now that $B\in \cB$ is a best approximation to $A$
in $\cB$. By the results stated above there is a $\psi \in \cFo$ such that
$\psi \in \cB^\perp$ and $\psi(A) = \|A-B\|$. Since $\psi$ is in the boundary
of the compact convex set $\cFo$, there exist extreme points,
$\psi_1, \dots, \psi_k$ of $\cFo$, with $k \leq p+1$, and there exist
positive real numbers $t_1, \dots, t_k$ with $\sum t_j = 1$
such that $\psi = \sum t_j\psi_j$. Each $\psi_j$ can be extended to
an extreme point of $\cAo$. We denote these extensions again by
$\psi_j$. Then we extend $\psi$ to $\cA$ by setting 
$\psi = \sum t_j\psi_j$. Note that $\|\psi\| = 1$, $\psi \in \cB^\perp$,
and $\psi(A) = \|A-B\|$. Since $|\psi_j(A-B)| \leq \|A-B\|$ for each $j$,
we must have $\psi_j(A-B) =  \|A-B\|$ for each $j$. (But we can not
expect that $\psi_j \in \cB^\perp$.) In this way we obtain:

\begin{proposition}
\label{prop-ext}
Let $\cA$ and $\cB$ be as above, with $p = \dim(\cB)$. Let $A \in \cA$
and $B \in \cB$ be given. Then $B$ is a best approximation to $A$ in
$\cB$ if and only if there are $k$ extreme points, 
$\psi_1, \dots, \psi_k$, of $\cAo$, with $k \leq p+1$, and there are
positive real numbers $t_1, \dots, t_k$ with $\sum t_j = 1$,
such that when we set $\psi = \sum t_j\psi_j$, we have $\psi \in \cAo$,
$\psi \in \cB^\perp$, and $\psi(A) = \|A-B\|$. Furthermore, 
$\psi_j(A) = \|A-B\|$ for each $j$. 
\end{proposition}


\section{Witnesses for best approximation in C*-algebras}
\label{wit-alg}

In this section we assume for simplicity that $\cA$ is a unital 
C*-algebra and that
$\cB$ is a unital C*-subalgebra of $\cA$ (with $1_\cA \in \cB$).
We let $\cAH$ denote the subspace of Hermitian elements of $\cA$,
and similarly for $\cBH$. It is easily seen that if $A \in \cAH$ and $A$
has a best approximation $B$ in $\cB$, then it has a best approximation
in $\cB$ that is actually in $\cBH$ and of no greater norm.
If $A \neq A^*$, then we can consider the Hermitian element
$
\begin{pmatrix}  0 & A^*  \\
                          A & 0      \end{pmatrix} 
$
in $M_2(\cA)$. If this element has a best approximation in $M_2(\cB)$,
then it is easily seen that it has a best approximation in $M_2(\cB)$
that is of the form 
$
\begin{pmatrix}  0 & B^*  \\
                          B & 0      \end{pmatrix} 
$
for $B \in \cB$ of no greater norm, so that $B$ is a best approximation
to $A$ in $\cB$ of no greater norm. In this way we can reduce the study
of best approximation in C*-algebras to that for Hermitian elements.
One does not need to take this path  ---  one can, for example, instead
work directly with elements of $\cAo$, much as we did
in Section \ref{secC^*}. But the path via Hermitian
elements involves somewhat more familiar arguments, using states, as 
we will see. Note that $\cAH$ is a vector space over $\bR$. 

We now
apply the results of the previous section. We seek witnesses for
elements of $A \in \cAH$ to be $\cB$-minimal, since if $A \in \cAH$
and $B \in \cBH$ then $B$ is a best approximation to $A$ in $\cB$
exactly if $A-B$ is $\cB$-minimal. So suppose now that $Z \in \cAH$
(with $Z \neq 0$), and that $Z$ is $\cB$-minimal. 
By rescaling, we see that it is
sufficient to treat the case in which $\|Z\| = 1$ We now assume this,
as it slightly simplifies the bookkeeping. Then, as seen in the previous
section, there exists a $\psi \in \cAD$ such that $\psi \in (\cBH)^\perp$
and $\psi(Z) = 1$. Note that necessarily $\|\psi\| = 1$. We can 
extend $\psi$ (uniquely) to be
a Hermitian element of $\cA'$ (of same norm). 

We are now exactly at the point where we can take advantage of
our C*-algebra setting. Let $\psi = (\psi^+ - \psi^-)/2$ be the Jordan
decomposition of $\psi$, as discussed, for example, in theorem 4.3.6 and 
remark 4.3.12 of \cite{KR1}. Thus $\psi^+$ and $\psi^-$ are
positive linear functionals on $\cA$ that are orthogonal in the sense
that 
\[
\|\psi^+\| \ + \ \|\psi^-\| \ = \ \|\psi^+ \ - \ \psi^-\| \quad (= 2\|\psi\| = 2)   .
\]
Now by assumption $1_\cA \in \cB$, and so
\[
0 = \psi(1_\cA) = (\psi^+(1_\cA) \ - \ \psi^-(1_\cA))/2   ,
\]
so that $\psi^+(1_\cA) = \psi^-(1_\cA)$, and consequently
$\|\psi^+\| = \|\psi^-\|$.
It follows that each of $\psi^+$ and $\psi^-$
is a state of $\cA$. Now
\[
2 = 2\psi(Z) = \psi^+(Z) - \psi^-(Z)   .
\]
Since $\|Z\| = 1$, it follows that $\psi^+(Z) = 1$ while $\psi^-(Z) = -1$.
This is an echo of the Chebyshev equi-oscillation phenomenon
mentioned earlier, with more echo to come in the next section. Because
$\psi^-(Z) = -1$ we have
\[
\psi^-((Z + 1_\cA)^2) = \psi^-(Z^2) + 2\psi^-(Z) + 1 = \psi^-(Z^2) - 1   .
\]
Since the left-hand side is clearly non-negative and the right-hand
side is clearly non-positive, we see that 
\[
\psi^-((Z + 1_\cA)^2) = 0 \quad \mathrm{and} \quad 
\psi^-(Z^2) = +1 = (\psi^-(Z))^2  .
\]
A similar calculation shows that 
\[
\psi^+((Z - 1_\cA)^2) = 0 \quad \mathrm{and} \quad
\psi^+(Z^2) = +1 = (\psi^+(Z))^2  .
\]
This means that each of $\psi^+$ and $\psi^-$ is ``definite'' on $Z$, as defined,
for example, in exercise 4.6.16 of \cite{KR1} (and see also 4.6.16
of \cite{KR3}). I thank Dick Kadison for pointing out to me the relevance of this 
definition to the present situation. Another way of expressing this
definiteness is that the ``mean-square deviation'' of $Z$ for 
$\psi^+$ and for $\psi^-$, as considered in remark 2.2.23(4) of \cite{Thr3}, is $0$.
As follows from exercise 4.6.16 of \cite{KR1}, this definiteness implies that
\[
\psi^+(ZC) = \psi^+(C) \quad \mathrm{and} \quad \psi^-(ZC) = - \psi^-(C)
\]
for all $C \in \cA$, as can be seen by applying the Cauchy-Schwarz
inequality to $\psi^+((Z - 1_\cA)C)$, and similarly for $\psi^-$. 

Let us now set $\var = (\psi^+ + \psi^-)/2$, which is a state on $\cA$.
In analogy with measure theory, we could call $\var$ the ``total
variation'' of $\psi$ and write $\var = |\psi|$.
From the calculations done some lines above we see that $\var(Z^2) = 1$
(from which one can show that $\var$ is definite on $Z^2$). But
we also see from these calculations that for every $B \in \cBH$ we have
\[
2\var(ZB) = \psi^+(ZB) + \psi^-(ZB) = \psi^+(B) - \psi^-(B)
= \psi(B) = 0,
\]
so that $\var(ZB) = 0$ for all $B \in \cBH$. Since $Z$ is Hermitian,
it follows that also $\var(BZ) = 0$ for all $B \in \cBH$, and so\[
\var(ZB + B^*Z) = 0
\]
for all $B \in \cB$. Note that $ZB + B^*Z$ is Hermitian. We arrive
in this way at one direction of the proof of the next theorem, which is 
essentially theorem 2.2 of \cite{AMM}, which in turn has antecedents
in section 5 of \cite{DMR}. The authors of these two papers have a purpose 
quite different from ours, namely to understand geodesics in the
homogeneous spaces $\cU_\cA/\cU_\cB$ where $\cU_\cA$ denotes
the unitary group of the C*-algebra $\cA$, and similarly for $\cU_\cB$.
They use the Finsler metric from the operator norm.The authors
give a direct proof of the next theorem, a proof that is somewhat more
complicated than what we have done above, and they provide almost no
indication that this theorem has anything to do with best
approximation, much less that it relates strongly to basic central 
results in the literature on best approximation along the lines that
we have discussed above. The authors do examine some interesting
specific finite-dimensional examples in which $\cA$ is a full matrix algebra
and $\cB$ is the subalgebra of diagonal matrices, or a subalgebra of
block-diagonal matrices.

\begin{theorem}
\label{thm-wit}
Let $\cA$ be a unital C*-algebra and let $\cB$ be a unital C*-subalgebra
of $\cA$. Let $A \in \cAH$. If $A$ is $\cB$-minimal, then there exists 
a state $\var$ of $\cA$ such that $\var(A^2) = \|A\|^2$ and
$\var(AB +B^*A) = 0$ for all $B \in \cB$. Conversely, if $A \in \cAH$ and
if there is a state $\var$ of $\cA$ such that  $\var(A^2) = \|A\|^2$ and
$\var$ satisfies the slightly weaker condition that $\var(AB +BA) = 0$ 
for all $B \in \cBH$, then $A$ is $\cB$-minimal.
\end{theorem}

\begin{proof}
By scaling we can assume that $\|A\| = 1$. We have given above the proof 
that if $A$ is $\cB$-minimal then there exists a state $\var$ satisfying the given 
conditions. Conversely, suppose that there exists a state $\var$ that satisfies 
the slightly weaker conditions. Then for any $B \in \cBH$ we have
\begin{align*} 
\var((A-B)^*(A-B)) & = \var(A^2) -\var(BA+AB)+\var(B^2)   \\
   & = \|A\|^2 + \var(B^2) \geq \|A\|^2  .
\end{align*}
Consequently $\|A-B\|^2 \geq \|A\|^2$, so that $A$ is $\cB$-minimal.
\end{proof}

Thus a state satisfying the (weaker) conditions of this theorem can
serve as a witness for the $\cB$-minimality of $A$. Since $\var$ is a state
and $A^2$ is positive, this may be more convenient to use than
our earlier $\psi$.

We remark that the condition that $\var(AB +BA) = 0$ 
for all $B \in \cBH$ is natural from the following point of view. Let
$\var$ be any state of $\cA$, let
$B \in \cB$ be given, and define a function, f, on $\bR$ by
\[
f(t) = \var((A+tB)^2)   .
\]
If $A$ is $\cB$-minimal, we can hope that $f$ takes a minimum value
at $t = 0$. But 
\[
f'(t) = \var(AB + BA) + 2t\var(B^2)   ,
\]
and so if $f'(0) = 0$ then $\var(AB+BA)=0$. Notice that 
$f''(0) = 2\var(B^2) \geq 0$. But these calculations do not
provide part of a proof for the theorem.

The proof of Theorem \ref{thm-wit} can be used to give some
information about the uniqueness of best approximations:

\begin{theorem}
\label{thm-uniq}
Let $\cA$ be a unital C*-algebra and let $\cB$ be a unital C*-subalgebra
of $\cA$. Let $A \in \cAH$. Suppose that $\var$ is a state of $\cA$ 
such that  $\var(A^2) = \|A\|^2$, and that
$\var$ satisfies the condition that $\var(AB +BA) = 0$ 
for all $B \in \cBH$, so that $A$ is $\cB$-minimal.
If the restriction of $\var$ to $\cB$ is faithful, then
0 is the unique best approximation to $A$ in $\cB$.
\end{theorem}

\begin{proof}
The last display of the proof of Theorem \ref{thm-wit} gives
\[
\var((A-B)^*(A-B)) = \|A\|^2 + \var(B^2)  .
\]
Consequently if the restriction of $\var$ to $\cB$ is faithful, 
then for $B \in \cBH$ with $B \neq 0$ we have $\var(B^2) > 0$,
so that $\|A-B\| > \|A\|$.
\end{proof}

\begin{corollary}
\label{corbest}
Let $\cA$ be a unital C*-algebra and let $\cB$ be a unital C*-subalgebra
of $\cA$. Let $A \in \cAH$. Suppose that $A$ has a best
approximation, $B$, in $\cBH$. Let $\var$ be a state of $\cA$ 
such that  $\var((A-B)^2) = \|A-B\|^2$ and
$\var((A-B)D +D(A-B)) = 0$ 
for all $D \in \cBH$. (Such a state is guaranteed 
to exist by Theorem \ref{thm-wit}.)
If the restriction of $\var$ to $\cB$ is faithful, then
$B$ is the unique best approximation to $A$ in $\cB$.
\end{corollary}

Thus in the case when $\var$ is not unique, it can be useful to
seek a $\var$ that is faithful, or at least has support as large as possible. 

Recall that we are interested in upper bounds on
the norms of best approximations.
Knowing $\var(B^2)$ gives us a lower bound on $\|B\|$, but that
may not be very useful.

For the present special context we can give a considerably simpler
proof of the key part of Theorem \ref{thmArv}:

\begin{corollary}
Let $\cA$ be a unital C*-algebra and let $\cB$ be a unital C*-subalgebra
of $\cA$. Let $A \in \cAH$. Suppose that $A$ has a best
approximation, $B$, in $\cBH$. Then there is a non-degenerate $*$-representation,
$(\cH, \pi)$, of $\cA$, and a Hermitian unitary operator $U \in \cL(\cH)$,
such that 
 \[
 L(A) = (1/2)\|[U, A]\|
 \quad \mathrm{while} \quad (1/2)\|[C, U]\| \leq L(C)  
\]
for all $C \in \cA$.
\end{corollary}

\begin{proof}
Let $Z = A - B$, so that $Z$ is $\cB$-minimal. Let $\var$ be the state
whose existence is guaranteed by Theorem \ref{thm-wit},  and let
$(\cH, \pi, \xi)$ be the GNS representation for $\var$.  Let $P$
be the orthogonal projection of $\cH$ onto the closure of
$\pi(\cB)\xi$. As before, it is clear that we have $[P, \pi(B)] = 0$
 for all $B \in \cB$. Notice that $\pi(Z)\xi \perp \pi(\cB)\xi$
 because $Z$ is Hermitian and $\var(ZB) = 0$ for all $B \in \cB$. Thus
 $P(Z\xi) = 0$. Also, because $\cB$ is unital we see that $\xi$
 is in $\pi(\cB)\xi$, so that $P(\xi) = \xi$. Finally, $\|Z\xi\| = \|Z\|$
 because $\var(Z^2) = \|Z\|^2$. Set $\eta = \pi(Z)\xi$, and
 let $U = 2P - I$. Then by the same calculations as done in the
 proof of Theorem \ref{thmArv} we find that 
 \[
 L(A) = (1/2)\|[U, A]\|
 \quad \mathrm{while} \quad (1/2)\|[C, U]\| \leq L(C)  
\]
for all $C \in \cA$.   
\end{proof}


\section{The case when $\cB$ is finite-dimensional}
\label{fin-dem}

In this section we assume that $\cA$ is a unital C*-algebra and that $\cB$ 
is a unital C*-subalgebra that is finite-dimensional. In this case best 
approximations always exist. We let $p = \dim(\cB)$. (We remark that
if, instead, $\cB$ is a closed two-sided ideal of $\cA$, then again
best approximations always exist. This is an immediate consequence
of proposition II.5.1.5 of \cite{Blk2}. But I have not seen how to
make use of this fact in the present context.)

The case in which $\cA$ itself is finite-dimensional is already quite
interesting. In particular, for any natural numbers $m$ and $n$ we
can set $\cA = (M_n(\bC))^m$, the C*-algebra of $m$-tuples of
 elements of $M_n(\bC)$. Equally well, $\cA$ can be viewed 
 as $M_n(\cD)$ where $\cD$ is the commutative C*-algebra of 
 functions on a set with $m$ elements. Then, along the lines
 of our comments in the introduction concerning matricial
 norms, it is natural to take $\cB$ to be the C*-subalgebra of
constant $m$-tuples, so that $\cB \cong M_n(\bC)$.
This approximation problem 
can be viewed as the problem of finding a center of a smallest
ball in $M_n(\bC)$ containing the $m$ entries of a
given element of $\cA$.
One might try taking the average of the $m$ elements, or at least
expect that the center will be in the convex hull of the $m$ elements.
But there is a striking theorem of Garvaki \cite{Grv} that says that if 
$\cV$ is a Banach space of dimension at least 3, and if for every three 
points of $\cV$ there exists a center for a ball of smallest radius 
containing the three points, such that this center lies in a plane 
containing the three point, then $\cV$ is a Hilbert space. Thus in
our C*-algebra setting we should not expect that centers will even 
be in the convex hull of the points. Interested readers can try their
hand at finding a center (that is, a best approximation in $\cB$)
for the case in which $m=3$, $n=2$, and
\[
A = \left\{ \begin{pmatrix}  1 & 0 \\
                                          0 & 0       \end{pmatrix}     \ , \ 
               \begin{pmatrix}  0 & 0  \\
                                         0 & 1             \end{pmatrix}     \ , \
              \begin{pmatrix}   0 & 1  \\                                              
                                         1 & 0             \end{pmatrix}     
     \right\}  .         
\]     
This example is small enough that it is not difficult, but it does seem
to illustrate fairly well some of the challenges of this kind of problem.

Let us now apply some of the results of Section \ref{wit-Ban} to
refine those of Section \ref{wit-alg}. Let $Z$ be an element of
$\cAH$ that is $\cB$-minimal.  We assume
as before that $\|Z\| = 1$. Much as in the discussion leading to Proposition
\ref{prop-ext}, we set $\cF = \cBH \oplus \bR Z$. Then, as seen in that
discussion, there exist extreme points
$\psi_1, \dots, \psi_k$ of $\cFo$ with $k \leq p+1$, and there exist
positive real numbers $t_1, \dots, t_k$ with $\sum t_j = 1$,
such that when we set $\psi = \sum t_j\psi_j$ we have $\psi(Z) = \|Z\| = 1$
and $\psi \in (\cBH)^\perp$. 
Each $\psi_j$ can be extended to be
an extreme point of $(\cAH)'_1$ and extended further to be a
Hermitian extreme point of $\cAo$. We denote these extensions again by
$\psi_j$. Then we extend $\psi$ to $\cA$ by setting 
$\psi = \sum t_j\psi_j$. Note that $\|\psi\| = 1$,  $\psi(Z) = 1$, and 
$\psi \in \cB^\perp$. 

Of course
\[
1 = \|Z\| = \psi(Z) =  \sum t_j\psi_j(Z)   ,
\]
and since $|\psi_j(Z)| \leq 1$ for each $j$, 
it follows that $\psi_j(Z) = 1$ for each $j$. But when we form the Jordan decompositions of the $\psi_j$'s, theorem 4.3.6 of \cite{KR1}, much as we
did in the previous section, the fact that each $\psi_j$ is an extreme point
tells us that either $\psi_j$ is a pure state of $\cA$ or the negative of
a pure state. Thus for each $j$ there is an integer $\e_j$ which is either $+1$
or $-1$ such that $\e_j\psi_j$ is a pure state. We set $\var_j = \e_j\psi_j$
for each $j$. Thus each $\var_j$ is a pure state, and $\var_j(Z) = \e_j$
for each $j$. This is again an echo of the Chebyshev equi-oscillation 
phenomenon mentioned
earlier, with a bit more echo to come.
Of course $\psi = \sum t_j\e_j\var_j$. 

Because $\var_j(Z) = \e_j$ for each
$j$, we can argue exactly as in the previous section to conclude that
for each $j$
\[
\var_j(Z^2) = 1 = \|\var_j(Z)\|^2  ,
\]
and that consequently $\var_j$ is ``definite'' on $Z$, so that
$\var_j(ZC) = \e_j\var_j(C)$ for every $C \in \cA$.

Let us now set $\var = \sum t_j\var_j$, so that $\var$ is a state of
$\cA$, not usually pure. Then for every $B \in \cBH$ we have
\begin{align*}
\var(ZB) &= \sum t_j\var_j(ZB) = \sum t_j\e_j\var_j(B)  \\
& = \sum t_j\psi_j(B) = \psi(B) = 0  .
\end{align*}

Since $Z$ is Hermitian, this implies, as in the previous section, that
$\var(ZB + B^*Z) = 0$ for all $B \in \cB$. We are now exactly in the
situation for one direction of Theorem \ref{thm-wit}, except that we now have 
the additional information about the decomposition of $\var$ into 
pure states that are definite on $Z$.

It is natural to set 
\[
\psi^+ = 2\sum\{t_j\psi_j: \e_j = +1\} \quad \mathrm{and} \quad 
\psi^- = 2\sum\{t_j\psi_j: \e_j = -1\} ,
\]
for then $\psi = (\psi^+ - \psi^-)/2$, while $\psi^+$ and $\psi^-$
are positive and $\var = (\psi^+ + \psi^-)/2$. Since $1_\cA \in \cB$
so that $\psi(1_\cA) = 0$, we see that $\psi^+(1_\cA) = \psi^-(1_\cA)$,
so that $\psi^+$ and $\psi^-$ must be states. Since $\|\psi\| = 1$, we see
that $\psi^+$ and $\psi^-$ are orthogonal, and so give the Jordan
decomposition of 2$\psi$. It follows easily that $\psi^+(Z) = 1$
while $\psi^-(Z) = -1$, so that each of $\psi^+$ and $\psi^-$
is definite on $Z$, as is to be expected.

We summarize part of the above discussion as follows:

\begin{theorem}
\label{thm-fd}
Let $\cA$ be a unital C*-algebra and let $\cB$ be a unital C*-subalgebra
of $\cA$ that is finite-dimensional, of dimension $p$. Let $Z$ be an
element of $\cAH$ that is $\cB$-minimal. Then
there exist pure states $\var_1, \dots, \var_k$ of $\cAo$
with $k \leq p+1$, such that each $\var_j$ is definite on $Z$, and there are
positive real numbers $t_1, \dots, t_k$ with $\sum t_j = 1$,
such that when we set $\var = \sum t_j\var_j$ then
$\var(Z^2) = \|Z\|^2$ and $\var(ZB + B^*Z) = 0$
for all $B \in \cB$.
\end{theorem}

Thus we are assured that we can always find a witness, $\var$, for 
the minimality of $Z$ that is expressed in terms of pure states definite on
$Z$ in the way stated in the theorem. This suggests that one way to
find a witness $\var$ for the minimality of $Z$
is to examine the pure states that are definite
on $Z$. We use this approach in the next section.


\section{The failure of the same-norm property}
\label{ex}

\begin{example}
\label{badnear}
As earlier, we let $\cA = (M_2(\bC))^3$, and we let $\cB$ be its subalgebra
of constant 3-tuples. We will actually just work with real Hermitian matrices, 
and it is easily seen that if $A \in \cAH$ and if its 3 entries are all real matrices,
then there will be a best approximation to $A$ in $\cB$ that is a real 
symmetric matrix.
We let
\[
Z = \left\{ \begin{pmatrix}  2 & 0 \\
                                          0 & 5       \end{pmatrix}     \ , \ 
               \begin{pmatrix}  4 & -3  \\
                                         -3 & -4             \end{pmatrix}     \ , \
              \begin{pmatrix}   4 & 3  \\                                              
                                         3 & -4             \end{pmatrix}     
     \right\}  .         
\]     
We do not require that $\|Z\| = 1$ so that we can work with integer entries. Then
\[
Z^2 = \left\{ \begin{pmatrix}  4 & 0 \\
                                          0 & 25       \end{pmatrix}     \ , \ 
               \begin{pmatrix}  25 & 0  \\
                                         0 & 25             \end{pmatrix}     \ , \
              \begin{pmatrix}   25 & 0  \\                                              
                                         0 & 25             \end{pmatrix}     
     \right\}  .         
\]     
We now use the tools developed in the previous section to prove that $Z$
is minimal. There are 5 pure states on $\cA$ that are definite on $Z$
and take value 25 on $Z^2$. Each of the entries of $Z$ is self-adjoint, 
and we use the pure states corresponding to the eigenvectors of the
entries for the eigenvalues $\pm 5$ of $Z$. We denote these states by $\var^j_\e$
where $j$ corresponds to the entry used, so $j = 1, 2, 3$, and $\e$ 
is $+$ or $-$ depending on the sign of the corresponding eigenvalue. 
For instance, $\var^1_+$ is the state that acts on the first entry of the elements
of $\cA$ and is determined by the eigenvector 
$(\begin{smallmatrix}  0  \\  1  \end{smallmatrix} )$
of the first entry of $Z$, while
$\var^2_-$ is the state that acts on the second entry of elements of $\cA$
and is determined by the eigenvector for the eigenvalue $-5$ of the
second entry of $Z$. Thus $\var^2_-(Z) = -5$. 

We seek $\psi$ and $\var$ with $\var = |\psi|$, satisfying the conditions
of Theorem \ref{thm-fd}  , expressed as a convex combination of the 5
pure states. Thus we seek coefficients that must be non-negative and sum to
1. The subspace of $\cBH$ consisting of real matrices has dimension 3,
and so we expect to need only 4 of the 5 pure states. This indicates
that $\psi$ and $\var$, if they exist, may not be unique, and calculations
show that this is indeed the case. In particular, calculations show that 
there do exist solutions, and that one of them is
\[
\psi = (1/18)(8\var^1_+ + \var^2_+ - 4 \var^2_- - 5\var^3_-)  ,
\]
so that $\psi^+ = (1/9)(8\var^1_+ + \var^2_+)$ and
$\psi^- = (1/9)(4 \var^2_- + 5\var^3_-)$  (for our conventions),
and $\var = (\psi^+ + \psi^-)/2$.
Once found, it is not hard to check that this is in fact a solution.
Thus $Z$ is minimal. Furthermore, $\var$ 
restricted to $\cB$ is faithful, because here this amounts to the
fact that the set of eigenvectors defining the pure states used above
to express $\psi$ spans $\bR^2$.  Then Theorem \ref{thm-uniq}
tells us that 0 is the unique nearest element
to $Z$ in $\cB$.

Now let $B =  \begin{pmatrix}   -8 & 0  \\                                              
                                         0 & 0             \end{pmatrix} $ , viewed
                                         as a constant 3-tuple in $\cBH$. Let $A = Z + B$.
It is easily seen that $\|A\| = 7$. But by the uniqueness of the
best approximation to $Z$, the best approximation to $A$ in $\cB$ is
unique, and is just $B$. But $\|B\| = 8$. Thus we see that $A$ has
no best approximation of norm no bigger than $\|A\|$. We see in
this way the failure of the same-norm approximation property.                                        
\end{example}                                   
                                        



\def\dbar{\leavevmode\hbox to 0pt{\hskip.2ex \accent"16\hss}d}
\providecommand{\bysame}{\leavevmode\hbox to3em{\hrulefill}\thinspace}
\providecommand{\MR}{\relax\ifhmode\unskip\space\fi MR }
\providecommand{\MRhref}[2]{%
  \href{http://www.ams.org/mathscinet-getitem?mr=#1}{#2}
}
\providecommand{\href}[2]{#2}



\begin{thebibliography}{99}


\bibitem{AMM}
Esteban Andruchow, Luis~E. Mata-Lorenzo, Alberto Mendoza, L{\'a}zaro Recht, and
  Alejandro Varela, \emph{Minimal matrices and the corresponding minimal curves
  on flag manifolds in low dimension}, Linear Algebra Appl. \textbf{430}
  (2009), no.~8-9, 1906--1928. \MR{2503942 (2010b:46133)}

\bibitem{Arv}
William Arveson, \emph{Interpolation problems in nest algebras}, J. Functional
  Analysis \textbf{20} (1975), no.~3, 208--233. \MR{0383098 (52 \#3979)}

\bibitem{BhS}
Rajendra Bhatia and Peter {\v{S}}emrl, \emph{Orthogonality of matrices and some
  distance problems}, Linear Algebra Appl. \textbf{287} (1999), no.~1-3,
  77--85, Special issue celebrating the 60th birthday of Ludwig Elsner.
  \MR{1662861 (99k:15042)}

\bibitem{Blk2}
B.~Blackadar, \emph{Operator algebras}, Encyclopaedia of Mathematical Sciences,
  vol. 122, Springer-Verlag, Berlin, 2006, Theory of $C{^{*}}$-algebras and von
  Neumann algebras, Operator Algebras and Non-commutative Geometry, III.
  \MR{2188261 (2006k:46082)}

\bibitem{ChL}
Man-Duen Choi and Chi-Kwong Li, \emph{The ultimate estimate of the upper norm
  bound for the summation of operators}, J. Funct. Anal. \textbf{232} (2006),
  no.~2, 455--476. \MR{2200742 (2006j:47010)}

\bibitem{Chr}
Erik Christensen, \emph{Perturbations of operator algebras. {II}}, Indiana
  Univ. Math. J. \textbf{26} (1977), no.~5, 891--904. \MR{0512368 (58
  \#23628b)}

\bibitem{Dvs}
Philip~J. Davis, \emph{Interpolation and approximation}, Blaisdell Publishing
  Co. Ginn and Co. New York-Toronto-London, 1963. \MR{0157156 (28 \#393)}

\bibitem{Dgl}
Ronald~G. Douglas, \emph{Banach algebra techniques in operator theory}, second
  ed., Graduate Texts in Mathematics, vol. 179, Springer-Verlag, New York,
  1998. \MR{1634900 (99c:47001)}

\bibitem{DMR}
Carlos~E. Dur{\'a}n, Luis~E. Mata-Lorenzo, and L{\'a}zaro Recht, \emph{Metric
  geometry in homogeneous spaces of the unitary group of a {$C^\ast$}-algebra.
  {II}. {G}eodesics joining fixed endpoints}, Integral Equations Operator
  Theory \textbf{53} (2005), no.~1, 33--50. \MR{2183595 (2007a:58006)}

\bibitem{Gjn}
P.~Gajendragadkar, \emph{Norm of a derivation on a von {N}eumann algebra},
  Trans. Amer. Math. Soc. \textbf{170} (1972), 165--170. \MR{0305090 (46
  \#4220)}

\bibitem{Grv}
A.~L. Garkavi, \emph{On the \v {C}eby\v sev center and convex hull of a set},
  Uspehi Mat. Nauk \textbf{19} (1964), no.~6 (120), 139--145. \MR{0175035 (30
  \#5221)}

\bibitem{HLT}
D.~Hadwin, D.~R. Larson, and D.~Timotin, \emph{Approximation theory and matrix
  completions}, Linear Algebra Appl. \textbf{377} (2004), 165--179.
  \MR{2021609 (2004j:47028)}

\bibitem{KR1}
R.~V. Kadison and J.~R. Ringrose, \emph{Fundamentals of the theory of operator
  algebras. {V}ol. {I}}, American Mathematical Society, Providence, RI, 1997,
  Elementary theory, Reprint of the 1983 original. \MR{98f:46001a}

\bibitem{KR3}
\bysame, \emph{Fundamentals of the theory of
  operator algebras. {V}ol. {III}}, Birkh\"auser Boston Inc., Boston, MA, 1991,
  Special topics, Elementary theory---an exercise approach. \MR{1134132
  (92m:46084)}

\bibitem{KR2}
\bysame, \emph{Fundamentals of the theory of operator algebras. {V}ol. {II}},
  Graduate Studies in Mathematics, vol.~16, American Mathematical Society,
  Providence, RI, 1997, Advanced theory, Corrected reprint of the 1986
  original. \MR{1468230 (98f:46001b)}

\bibitem{Pd2}
Gert~Kjaerg{\.a}rd Pedersen, \emph{\v {C}eby\v sev subspaces of {$C^{\ast}
  $}-algebras}, Math. Scand. \textbf{45} (1979), no.~1, 147--156. \MR{567440
  (81i:46077)}
  
  \bibitem{R17}
Marc~A. Rieffel,   
\emph{Vector bundles and {G}romov-{H}ausdorff distance}, J. K-Theory
 \textbf{5} (2010), 39--103, arXiv:math.MG/0608266.

\bibitem{R21}
\bysame, 
\emph{Leibniz seminorms for ``{M}atrix algebras converge to
  the sphere''}, Quanta of Maths, 543--578, Clay Mathematics Proceedings,
  vol.~11, Amer. Math. Soc., Providence, R.I., 2011, arXiv:0707.3229.

\bibitem{Rbt}
A.~Guyan Robertson, \emph{Best approximation in von {N}eumann algebras}, Math.
  Proc. Cambridge Philos. Soc. \textbf{81} (1977), no.~2, 233--236.
  \MR{0473860 (57 \#13519)}

\bibitem{RbY}
A.~Guyan Robertson and David Yost, \emph{Chebyshev subspaces of operator
  algebras}, J. London Math. Soc. (2) \textbf{19} (1979), no.~3, 523--531.
  \MR{540068 (80j:46097)}

\bibitem{Rdn}
Walter Rudin, \emph{Functional analysis}, second ed., International Series in
  Pure and Applied Mathematics, McGraw-Hill Inc., New York, 1991. \MR{1157815
  (92k:46001)}

\bibitem{Sng}
Ivan Singer, \emph{Best approximation in normed linear spaces by elements of
  linear subspaces}, Die
  Grundlehren der mathematischen Wissenschaften, Band 171, Publishing House of
  the Academy of the Socialist Republic of Romania, Bucharest, 1970.
  \MR{0270044 (42 \#4937)}

\bibitem{Stm}
Joseph~G. Stampfli, \emph{The norm of a derivation}, Pacific J. Math.
  \textbf{33} (1970), 737--747. \MR{0265952 (42 \#861)}

\bibitem{Thr3}
Walter Thirring, \emph{A course in mathematical physics. {V}ol. 3},
  Springer-Verlag, New York, 1981, Quantum mechanics of atoms and molecules,
  Translated from the German by Evans M. Harrell, Lecture Notes in Physics,
  141. \MR{625662 (84m:81006)}

\bibitem{Wuw1}
Wei Wu, \emph{Non-commutative metrics on matrix state spaces}, J. Ramanujan
 Math. Soc. \textbf{20} (2005), no.~3, 215--254, arXiv:math.OA/0411475.
  \MR{2181130}

\bibitem{Wuw2}
\bysame, \emph{Non-commutative metric topology on matrix state space}, Proc. Amer. Math. Soc. \textbf{134} (2006), no.~2, 443--453, arXiv:math.OA/0410587.
  \MR{2176013 (2006f:46072)}

\bibitem{Wuw3}
\bysame, \emph{Quantized {G}romov-{H}ausdorff distance}, J. Funct. Anal.
  \textbf{238} (2006), no.~1, 58--98, arXiv:math.OA/0503344. \MR{2234123}

\end{thebibliography}

\end{document}